\documentclass{amsart}
\usepackage{enumerate}
\usepackage{amssymb}
\usepackage{amsmath}
\usepackage{mathrsfs}
\usepackage{amsthm}
\usepackage{mathtools}
\usepackage{multirow}
\usepackage{float}

\newtheorem{definition}{Definition}
\newtheorem{theorem}{Theorem}

\newtheorem{corollary}[theorem]{Corollary}
\newtheorem{remark}[theorem]{Remark}

\newtheorem{question}[theorem]{Question}
\newtheorem{example}[theorem]{Example}
\usepackage[colorlinks,linkcolor=blue]{hyperref}

\allowdisplaybreaks[4]

\newcommand{\Sc}{{\rm Sc}}
\newcommand{\Ric}{{\rm Ric}}
\newcommand{\Sph}{\mathbb{S}}

\begin{document}
\title{Gap phenomenon for scalar curvature}
\author{Yukai Sun}
\address{Key Laboratory of Pure and Applied Mathematics,
School of Mathematical Sciences, Peking University, Beijing, 100871, P. R. China
}
\email{sunyukai@math.pku.edu.cn}

\author{Changliang Wang}
\address{
School of Mathematical Sciences and Institute for Advanced Study, Key Laboratory of Intelligent Computing and Applications(Ministry of Education), Tongji University, Shanghai 200092, China}
\email{wangchl@tongji.edu.cn}
\date{}
\begin{abstract}
Inspired by Goette-Semmelmann \cite{GSSU2002}, we derive an estimate for the scalar curvature without a nonnegativity assumption on curvature operator. As an application, we show that, on an even dimensional closed manifold with nonzero Euler characteristic, any Riemannian metric $g$ is $\epsilon$-gap distance extremal for some $\epsilon \geq 0$. For manifolds with boundary, inspired by Lott \cite{JL2021}, we obtained a similar estimate for scalar curvature and mean curvature. We apply the estimate on certain Euclidean domains to study a Gromov's question in \cite{GM20233} concerning the extension problem of metric on the boundary to the interior.
\end{abstract}
\maketitle

\section{Introduction}
In \cite{GM2018}, Gromov proposed the following concept of the $\epsilon$-gap distance extremality for a metric. 
\begin{definition}[Gromov \cite{GM2018}, page 152]\label{def-epsilon-gap-extremality}
    A Riemannian metric $g_{0}$ on a manifold $M^{n}$ is $\epsilon$-gap ($\epsilon\geq 0$) distance extremal if there is no metric $g$ on $M$ such that $g\geq g_{0}$ and
    \[\operatorname{Sc}_{g}-\operatorname{Sc}_{g_{0}}>\epsilon.\]
\end{definition}
In particular, for $\epsilon=0$, a $0$-gap distance extremal metric is said to be distance extremal. This extremality property of $g_0$ is equivalent to that any metric $g \geq g_0$ with $\Sc_g \geq \Sc_{g_0}$ must have $\Sc_g = \Sc_{g_0}$. If this further implies $g$ is isometric to $g_0$, then $g_0$ is said to be distance rigid. Moreover, by assuming $|v \wedge w|_g \geq |v \wedge w|_{g_0}$ for all $v, w \in \Gamma(TM)$ instead of $g \geq g_0$ in the above, one can correspondingly define area extremality and area rigidity property of metric $g_{0}$.

Many interesting results have been established concerning the extremality and rigidity of metrics. One of them is the following Llarull’s theorem, which was anticipated by Gromov in \cite{GM1986} on page 113 and proved by Llarull in \cite{ML1998}.
\begin{theorem}[Llarull \cite{ML1998}]\label{Thm-LLarull}
    Let $(M^{n},g)$ be a closed spin Riemannian manifold and $(\mathbb{S}^{n},g_{\mathbb{S}^{n}})$ the standard sphere with sectional curvature $1$. If there exists an area nonincreasing map $f:M^{n}\to \mathbb{S}^{n}$ with nonzero degree, then either the scalar curvature $\operatorname{Sc}_{g}$ of $g$ on $M$ satisfies
    \[\inf_{x\in M}\operatorname{Sc}_{g}(x)\leq n(n-1),\]
    or $M\cong \mathbb{S}^{n}$ and $f$ is an isometry.
\end{theorem}
We recall that a map $f:(N,g_{N})\to(M,g_{M})$ is distance nonincreasing if
\[|f_{\ast}v|_{g_{M}}\leq |v|_{g_{N}}\]
for all $x\in N$ and $v\in T_{x}N$. Additionally, $f:(N,g_{N})\to(M,g_{M})$ is area non-increasing if 
\[|f_{\ast}v\wedge f_{\ast} w|_{g_{M}}\leq |v\wedge w|_{g_{N}}\]
for all $x\in N$ and $v,w\in T_{x}N$.
Cecchini-Wang-Xie-Zhu \cite{CWXZ2024} recently proved Llarull's theorem for non-spin closed manifolds of dimension four, provided that the map $f$ is distance nonincreasing. 

Kramer \cite{WK2000} proved that for a Riemannian submersion with totally geodesic fibers, if the total space has distance rigidity property then so is the base space.
Goette-Semmelmann \cite{GSSU2002} considered the case of the target space $\mathbb{S}^{n}$ being replaced with a certain type of compact symmetric space with nonnegative curvature operator, and proved the area extremality of such symmetric spaces, and also area rigidity provided $0 < \Ric_g < \frac{\Sc_g}{2}g$. Bettiol-Goodman \cite{BG2024} proved area extremality of closed simply-connected four-dimensional Riemannian manifold with nonnegative sectional curvature, provided that certain linear combination of curvature operator and Hodge star operator is nonnegative; and they also proved area rigidity for such 4-manifolds, provided $0 < \Ric_g < \frac{\Sc_g}{2}g$. Zhu \cite{zhu2020}   proved the area-rigidity of standard product metric on $\mathbb{S}^2 \times T^n$ with $n+2 \leq 7$. Assume that $(M^n, g)$, $n\leq 8$, is a closed oriented Riemannian manifold, and a map $f: (M, g) \to (\Sph^k \times T^{n-k}, g_{\Sph^k} + g_{T^{n-k}})$ is area $\epsilon$-contracting map and nonzero degree, where $k=2$ or $3$. Cecchini-Schick \cite{CT2021} proved that the minimum of $\Sc_g \leq 2k\epsilon$. Hao-Shi-Sun \cite{HSS2024} proved the distance rigidity of $S^3 \times T^{n-3}$ with $n \leq 7$. Hu-Liu-Shi \cite{HLS2023} proved Llarull's theorem in dimension $3$ by using the $\mu$-bubble method. The rigidity of scalar curvature on manifolds with boundaries or corners has also been studied in \cite{BCHW,BS2024,CW2023,LC2024,JL2021,RD2023,Sun2024,WX20232,WX2023,WXY2021}. For the case of noncompact manifolds, relevant studies can be found in \cite{GMLB1983,HSS2024,zhang2020}. 

In contrast to the aforementioned scalar curvature extremality and rigidity results for certain metrics, a generic metric on a closed manifold does not have scalar curvature extremality or rigidity properties. For example, the manifold $\mathbb{S}^{2}\times \Sigma$ with the standard product metric $g_{0}=g_{\mathbb{S}^2}+g_{\Sigma}$ for the high genus $l\geq 2$ surface $\Sigma$ with constant curvature $-c(c>0)$, the metric $g_{0}$ is not extremal since the metric $g=g_{\mathbb{S}^2}+\lambda g_{\Sigma}(\lambda \geq 1)$ satisfies $g\geq g_{0}$ and $\Sc_{g}\geq\Sc_{g_{0}}$. In this example, when the distance increases, the scalar curvature increases as well. So it is a natural problem to estimate how much the scalar curvature can increase, i.e. estimate the upper bound of the infimum of the scalar curvature, or equivalently investigate the $\epsilon$-gap extremality property defined in Definition \ref{def-epsilon-gap-extremality}. For this problem, we obtain the following $\epsilon$-gap distance extramality result.

\begin{theorem}\label{Thm-gap}
    Let $M^{2n}$ be a closed $2n$-dimensional manifold with nonzero Euler characteristic. Then a complete Riemannian metric $g_{0}$ on $M$ is the $\epsilon$-gap distance extremal for any $\epsilon \geq 2n(2n-1)(|\mathcal{R}_{\min}| - \mathcal{R}_{\min})$, where $\mathcal{R}_{\min}$ is defined as in (\ref{eqn-Rmin}) for the metric $g_{0}$.
\end{theorem}

Theorem \ref{Thm-gap} follows from the following upper bound estimate for the infimum of the scalar curvature, which is inspired by Goette-Semmelmann's work \cite{GSSU2002} on manifolds with nonnegative curvature operator. We carry out the calculations for the estimate without any restriction on the curvature operator.



\begin{theorem}\label{Thm-main}
  Consider two closed, complete, connected spin Riemannian manifolds $(M^{2n},g_{0})$ and $(N^{2n},g)$.  If there exists a smooth area nonincreasing map $f: N\to M$ with nonzero degree and $M$ has nonzero Euler characteristic $\chi(M)$. Then the scalar curvature $\operatorname{Sc}_{g}$ of $g$ satisfies
  \begin{equation}\label{ineqn-Sc}
    \inf_{x\in N}\left[\operatorname{Sc}_{g}(x)-\operatorname{Sc}_{g_{0}}(f(x))\right]\leq -2n(2n-1)\mathcal{R}_{\min}+2n(2n-1)|\mathcal{R}_{\min}|  
  \end{equation}
  where, for the curvature operator $\mathcal{R}$ of the metric $g_0$ on $M$,
  \begin{equation}\label{eqn-Rmin}
  \mathcal{R}_{\min} := \inf_{x\in M}\{\text{the minimal eigenvalue of }\mathcal{R}_{x}:\wedge^2T_{x}M\to \wedge^2T_{x}M\}.
  \end{equation}
  Moreover, if
  \[\operatorname{Sc}_{g}(x)-\operatorname{Sc}_{g_{0}}(f(x))\geq -2n(2n-1)\mathcal{R}_{\min}+2n(2n-1)|\mathcal{R}_{\min}|.\]
  Then $\mathcal{R}_{\min}\geq 0$ and $\operatorname{Sc}_{g}(x)=\operatorname{Sc}_{g_{0}}(f(x))$.
\end{theorem}
\begin{remark}
{\rm
Theorem \ref{Thm-main} also holds if we only assume that $TN \oplus f^* TM$ admits a spin structure. This is equivalent to 
\begin{equation}\label{eqn-stiefel-Whitney-relation}
    w_2 (TN) = f^* (w_2(TM)),
\end{equation}
where $w_{2}(TM)$ is the second Stierel-Whitney c1ass of $TM$.
In particular, if $N=M$ and $f$ is the identity map, then the equation (\ref{eqn-stiefel-Whitney-relation}) always holds. In this case, by applying Theorem \ref{Thm-main}, we obtain Theorem \ref{Thm-gap}.
}
\end{remark}

The basic idea of the proof of Theorem \ref{Thm-main} is using the index theory of the Dirac operator on a twisted spinor bundle. This method dates back to Gromov-Lawson \cite{GMLB1983} and has been demonstrated to be very powerful in Llarull \cite{ML1998}, Goette-Semmelmann \cite{GSSU2002} and others. One of the key ingredients in this argument is the Lichnerowicz-Weitzenb\"ock-Bochner (LWB, for short) formula for the twisted Dirac operator, in which the term involving the curvature operator of the target space needs to be estimated. In \cite{ML1998}, Llarull dealt with the case of the target space being a round sphere with a constant positive curvature operator. In \cite{GSSU2002}, Goette-Semmelmann introduced a trick to estimate the term involving a nonnegative curvature operator of the target space. In this paper, following Goette-Semmelmann's trick in \cite{GSSU2002}, we estimate the curvature term in LWB formula but without any restriction on the curvature of the target space.

\begin{example}
{\rm
    Let us consider
 the manifold $(\mathbb
{S}^{2}\times \Sigma,g_{0}=g_{\mathbb{S}^2}+g_{\Sigma})$, where $(\Sigma, g_{\Sigma})$ is a genus $l\geq 2$ surface with constant sectional curvature $-c (c>0)$. It is clear that $\mathbb
{S}^{2}\times \Sigma$ does not admit any metric with either nonnegative or nonpositive curvature operator. The metric $g=g_{\mathbb{S}^2}+\lambda g_{\Sigma}\geq g_{0}$ for $\lambda\geq 1$. The scalar curvature of $g$ satisfies
\begin{equation*}
  \inf_{x\in \mathbb{S}^2\times \Sigma}\operatorname{Sc}_{g}(x)=2-\frac{2c}{\lambda}\leq  \operatorname{Sc}_{g_{0}}-2n(2n-1)\mathcal{R}_{\min}+2n(2n-1)|\mathcal{R}_{\min}|=2+22c,  
\end{equation*}
for any $c>0$. Note that for any fixed $\lambda >1$, $\lim\limits_{c\to0} 2 - \frac{2c}{\lambda} = \lim\limits_{c\to0} 2+22c = 2$. This indicates that the estimate in (\ref{ineqn-Sc}) is optimal.
}
\end{example}

Now we summarize the results in Gromov-Lawson \cite{GMLB1983}, Goette-Semmelmann \cite{GSSU2002}, and Theorem \ref{Thm-main} in the following table.

\begin{table}[H]
    \centering
    \begin{tabular}{|p{1.6cm}|p{3.2cm}|p{3.2cm}|p{3cm}|}
    \hline
    \multicolumn{4}{|c|}{$f:(N^{2n},g)\to (M^{2n},g_{0})$, \, $\chi(M)\neq 0$, \, $M$ and $N$ are spin} \\
\hline
    \multicolumn{1}{|c|}{$g_{0}$}&with nonnegative curvature operator&no curvature operator condition& with nonpositive curvature operator  \\
    \hline
   \multicolumn{1}{|c|}{$f$}  &area nonincreasing $\operatorname{deg}f\neq 0$&area nonincreasing $\operatorname{deg}f\neq 0$& $\operatorname{deg}f\neq 0$\\
   \hline
   Conclusion &$\inf\limits_{x\in N}(\operatorname{Sc}_{g}-\operatorname{Sc}_{g_{0}})\leq 0$&$\inf\limits_{x\in N}(\operatorname{Sc}_{g}-\operatorname{Sc}_{g_{0}})\leq c_{0}$&$\inf\limits_{x\in N}\operatorname{Sc}_{g}\leq 0$\\
   \hline
\end{tabular}
    \vskip0.1in
    \caption{Here $c_{0}=-2n(2n-1)\mathcal{R}_{\min}+2n(2n-1)|\mathcal{R}_{\min}|$. The second column is the result in \cite{GSSU2002}; The third column is Theorem \ref{Thm-main}; The fourth column is the result in \cite{GMLB1983}, since manifolds with nonpositive curvature operator are enlargeable.}
    \label{tab:my_label}
\end{table}

For manifolds with boundary, following Lott  \cite{JL2021}, we have
\begin{theorem}\label{Thm-manifold-boundary}
Let $(N^{2n},g)$ and $(M^{2n},g_{0})$ be compact connected spin manifolds with boundary. Let $f:N\to M$ be a smooth map and let $\partial f:\partial N\to \partial M$ denote the restriction to the boundary of $f$. Assume that
\begin{enumerate}[$(1)$]
    \item $f$ is area nonincreasing and $\partial f$ is distance nonincreasing,
    \item  $M$ has nonzero Euler characteristic and $f$ have nonzero degree.
\end{enumerate}
Then either
    \begin{equation}\label{eqn-scal-bdry}
    \inf_{x\in N}\left[\operatorname{Sc}_{g}(x)-\operatorname{Sc}_{g_{0}}(f(x))\right]\leq -2n(2n-1)\mathcal{R}_{\min}+2n(2n-1)|\mathcal{R}_{\min}|
    \end{equation}
    or 
    \begin{equation}\label{eqn-mean-curv-bdry}
        \inf_{y\in \partial N}\left[H_{ N}(y)-H_{ M}(\partial f(y))\right]\leq -(2n-1)A^{M}_{\min}+(2n-1)|A^{M}_{\min}|,
    \end{equation} 
    where $H_{N}$ and $H_{M}$ are the mean curvature (the trace of the second fundamental form) of $\partial N$ and $\partial M$ in $N$ and $M$, respectively, with respect to the outer unit normal vector, and for the adjoint operator $\mathcal{A}$ of the second fundamental form $A^{M}$ on $\partial M$,  
    \[A^{M}_{\min} :=\inf_{x\in\partial M}\{\text{the minimal eigenvalue of }\mathcal{A}_{x}:T_{x}\partial M\to T_{x}\partial M\}.\]
    Furthermore, if 
    \[
    \operatorname{Sc}_{g}(x)-\operatorname{Sc}_{g_{0}}(f(x))\geq -2n(2n-1)\mathcal{R}_{\min}+2n(2n-1)|\mathcal{R}_{\min}|
    \]   
    and 
    \[
        H_{ N}(y)-H_{ M}(\partial f(y))\geq -(2n-1)A^{M}_{\min}+(2n-1)|A^{M}_{\min}|.
    \]
    Then 
    \[\operatorname{Sc}_{g}(x)=\operatorname{Sc}_{g_{0}}(f(x)), \forall x \in N, \quad \mathcal{R}_{\min}\geq 0,\]
    and
    \[
        H_{ N}(y)-H_{ M}(\partial f(y))= -(2n-1)A^{M}_{\min}+(2n-1)|A^{M}_{\min}|.
    \]
    \[
    \begin{cases}
        H_{ N}(y)=H_{ M}(\partial f(y)), \forall y \in \partial N, & \text{ if } A^{M}_{\min}= 0\\
        \partial f:\partial N\to\partial M \text{ is local isometric }, & \text{ if } A^{M}_{\min}\neq 0.
    \end{cases}
    \]
\end{theorem}
As an application of Theorem \ref{Thm-manifold-boundary}, by letting $N=M$ be a domain in $\mathbb{R}^{2n}$ and $f$ be the identity map, we have
\begin{corollary}\label{cor-Euclidean domain}
    Let $\Omega \subset \mathbb{R}^{2n}$ be a compact Euclidean connected domain with smooth boundary and $\chi(\Omega) \neq 0$, $g$ be a Riemannian metric on $\Omega$ such that $\operatorname{Sc}_{g}\geq 0$. Then we have
\begin{eqnarray*}
    \inf_{x \in \partial \Omega}[H_{g}(x)-\hat{c}H_{g_{E}}(x)]\leq \hat{c}\left[-(2n-1)A^{M}_{\min}+(2n-1)|A^{M}_{\min}|\right],
\end{eqnarray*}
here $g_{E}$ is the Euclidean metric and 
\[\hat{c}=\sup \left\{\frac{|v|_{g_{E}}}{|v|_{g}} \, \mid \, v \neq 0\in T_{x}\partial \Omega, \, x \in \partial \Omega \right\},\]
\end{corollary}

This is related to the following question raised by Gromov, concerning the extension problem of metric on boundary to interior with prescribed lower bounds of scalar curvature and mean curvature of the boundary. 
\begin{question}[Gromov \cite{GM20233}, page 234]\label{Que-Gromov}
{\rm
    Let $X$ be a compact $n$-dimensional manifold with boundary, and let $Y_{i}\subset X$ be the connected components of the boundary.
    \begin{enumerate}[(a)]
        \item Given numbers $\sigma$ and $c_{i}$, when does there exist a Riemannian metric $g$ on $X$, such that $\operatorname{Sc}(X)\geq\sigma$ and the mean curvatures of $Y_{i}$ are bounded from below by $c_{i}$,
        \[H_{g}(Y_{i})\geq c_{i}?\]
        \item Let all $Y_{i}$ be diffeomorphic to the sphere $\mathbb{S}^{n-1}$ and, besides $\sigma$ and $c_{i}$, we are given positive numbers $\kappa_{i}$.

        When does there exist a Riemannian metric $g$ on $X$, such that $\operatorname{Sc}(X)\geq \sigma$, the induced metrics $g|_{Y_{i}}$ has constant sectional curvatures $\kappa_{i}$ and 
        \[H_{g}(Y_{i})\geq c_{i}?\]
        \item Let now Riemannian metrics $g_{i}$ on $Y_{i}$ be given. When does there exist a Riemannian metric $g$ on $X$, such that $\operatorname{Sc}_{g}\geq \sigma$, $g|_{Y_{i}}= g_{i}$ and
        \[H_{g}(Y_{i})\geq c_{i}?\]
    \end{enumerate}
}
\end{question}

For part (c) in Question \ref{Que-Gromov}, if $X$ is an even-dimensional Euclidean compact connected domain $(\Omega, g_{E})$ with a smooth boundary and $\chi(\Omega)\neq 0$, then there exists no metric $g$ on $X$ such that $\operatorname{Sc}_{g}\geq 0$, $g|_{\partial \Omega}=g_{E}|_{\partial\Omega} $ and 
\[H_{g}(x)>H_{g_{E}}(x)+\left[-(2n-1)A^{\Omega}_{\min}+(2n-1)|A^{\Omega}_{\min}|\right], \forall x \in \partial \Omega.\]
If we assume that the boundary $\partial \Omega$ is strictly convex, i.e. $A^{\Omega}_{\min} >0$, then this non-existence conclusion also follows from Corollary 1.2 in \cite{JL2021} and Corollary 1.7 in \cite{CHZ2024}.
There are also many interesting results about the estimate of the integral of mean curvature on manifolds nonnegative scalar curvature and boundary that can be convexly or mean-convexly embedded into the Euclidean space by Eichmair-Miao-Wang \cite{EMW2012}, Shi-Tam \cite{ST2002}, Shi-Wang-Wei \cite{SWW2022} and Shi-Wang-Wei-Zhu \cite{SWWZ2021}. 

The paper is organized as follows. In section \ref{sec-a}, we prove Theorem \ref{Thm-main}. In section \ref{sec-b}, we prove Theorem \ref{Thm-manifold-boundary}.

{\em Acknowledgements}: The authors would like to express their gratitude to Professor Xianzhe Dai and Professor Yuguang Shi for their interest and conversation about this work. Yukai Sun is partially funded by the National Key R\&D Program of China Grant 2020YFA0712800. Changliang Wang is partially supported by the Fundamental Research Funds for the Central Universities and Shanghai Pilot Program for Basic Research.


\section{Proof of Theorem \ref{Thm-main}}\label{sec-a}
In this section, we prove Theorem \ref{Thm-main}, basically following the argument in Goette-Semmelmann \cite{GSSU2002} with some modification to remove the nonnegativity constraint on the curvature operator. 
\begin{proof}
   Let $Cl_{2n}$ denote the Clifford algebra of $\mathbb{R}^{2n}$ with its standard inner product. Then $\mathbb{C}l_{2n} := Cl_{2n}\otimes_{\mathbb{R}} \mathbb{C} \cong End(\mathbb{C}^{2^n})$, and so $\mathbb{C}l_{2n}$ has a unique irreducible complex representation. By restricting to ${\rm Spin}_{2n}\subset \mathbb{C}l_{2n}$, one obtains a unitary representation $\rho$ of ${\rm Spin}_{2n}$. For a spin manifold $M$, there exists a principal ${\rm Spin}$-bundle $P_{{\rm Spin}_{2n}}$, and the associated complex vector bundle
\[S(M)=P_{{\rm Spin}_{2n}}\times_{\rho}\mathbb{C}^{2^n}\]
is called the spinor bundle over $M$. 

Similarly, on the spin manifold $N$, there exits the spinor bundle $S(N)$. Consider the bundle $S:=S(N)\otimes E$, where $E := f^{\ast}S(M)$, with the canonical tensor product connection. It is a twisted spin bundle over $N$. By Lichnerowicz-Weitzenb\"ock-Bochner formula for the twisted Dirac operator $D_E$ on $S$, there is 
\begin{equation}\label{eqn-Dirac}
D_{E}^2=\nabla^{\ast}\nabla+\frac{\operatorname{Sc}_{g}}{4}+\frac{1}{2}\sum_{j,k}(e_{j}\cdot e_{k}\cdot\sigma)\otimes R^{E}_{e_{j},e_{k}}(v)  
\end{equation}
for $\sigma\otimes v\in S$, where $R^{E}$ denotes the curvature tensor of the connection on $f^{\ast}E(M)$ and $\{e_{i}\}$ is a orthonormal basis of $T_{x}N$ for $x \in N$. Let $\{\epsilon_{i}\}$ be a local orthonomal basis of $T_{f(x)}M$. We now handle the term $\frac{1}{2}\sum_{j,k}(e_{j}\cdot e_{k}\cdot\sigma)\otimes R^{E}_{e_{j},e_{k}}(v)$. Define the Clifford action of two form by
$e_{i}\wedge e_{j}\cdot =e_{i}\cdot e_{j}\cdot$. Let $\bar{\omega}_{i}$, $\omega_{k}$ be the orthonormal basis of $\wedge^2T_{x}N$ and $\wedge^2T_{f(x)}M$, respectively. Then
\begin{eqnarray*}
    &&\frac{1}{2}\sum_{i,j=1}^{2n}(e_{i}\cdot e_{j}\cdot \sigma)\otimes (R^{E}_{e_{i}e_{j}}v)\\
    &=&\frac{1}{8}\sum_{i,j,k,l}g_{0}(R^{E}_{f_{\ast}(e_{i})f_{\ast}(e_{j})}\epsilon_{k},\epsilon_{l})(e_{i}\cdot e_{j}\cdot \sigma)\otimes (\epsilon_{k}\cdot\epsilon_{l} \cdot v)\\
    &=&-\frac{1}{8}\sum_{i,j,k,l}g_{0}(\mathcal{R}(f_{\ast}(e_{i})\wedge f_{\ast}(e_{j})),\epsilon_{l}\wedge\epsilon_{k})(e_{i}\wedge e_{j}\cdot \sigma)\otimes (\epsilon_{l}\wedge\epsilon_{k} \cdot v)\\
    &=&-\frac{1}{2}\sum_{i,j}g_{0}(\mathcal{R}(f_{\ast}\bar{\omega}_{i}),\omega_{j})(\bar{\omega}_{i}\otimes \omega_{j}) \cdot \sigma\otimes v\\
    &=&-\frac{1}{2}\sum_{i,j}g_{0}((\mathcal{R}-\mathcal{R}_{\min})(f_{\ast}\bar{\omega}_{i}),\omega_{j})(\bar{\omega}_{i}\otimes \omega_{j}) \cdot \sigma\otimes v\\
    &&-\frac{1}{2}\sum_{i,j}\mathcal{R}_{\min}g_{0}(f_{\ast}\bar{\omega}_{i},\omega_{j})(\bar{\omega}_{i}\otimes \omega_{j}) \cdot \sigma\otimes v.
\end{eqnarray*}
By a direct computing, we have
\begin{equation}\label{eqn-1-5}
    \begin{cases}
        g_{0}((\mathcal{R}-\mathcal{R}_{\min})(\epsilon_{i}\wedge\epsilon_{j}),\epsilon_{k}\wedge\epsilon_{l}) &=   g_{0}((\mathcal{R}-\mathcal{R}_{\min})(\epsilon_{k}\wedge\epsilon_{l}),\epsilon_{i}\wedge\epsilon_{j})\\
    g_{0}((\mathcal{R}-\mathcal{R}_{\min})(\epsilon_{i}\wedge\epsilon_{j}),\epsilon_{k}\wedge\epsilon_{l}) &=  -g_{0}((\mathcal{R}-\mathcal{R}_{\min})(\epsilon_{j}\wedge\epsilon_{i}),\epsilon_{k}\wedge\epsilon_{l})\\
    g_{0}((\mathcal{R}-\mathcal{R}_{\min})(\epsilon_{i}\wedge\epsilon_{j}),\epsilon_{k}\wedge\epsilon_{l}) &=   -g_{0}((\mathcal{R}-\mathcal{R}_{\min})(\epsilon_{i}\wedge\epsilon_{j}),\epsilon_{l}\wedge\epsilon_{k})\\
    g_{0}((\mathcal{R}-\mathcal{R}_{\min})(\epsilon_{i}\wedge\epsilon_{j}),\epsilon_{k}\wedge\epsilon_{l})
    &=-g_{0}((\mathcal{R}-\mathcal{R}_{\min})(\epsilon_{j}\wedge\epsilon_{k}),\epsilon_{i}\wedge\epsilon_{l})\\
    &\quad-g_{0}((\mathcal{R}-\mathcal{R}_{\min})(\epsilon_{k}\wedge\epsilon_{i}),\epsilon_{j}\wedge\epsilon_{l}).
    \end{cases}
\end{equation}

We claim that
\begin{eqnarray}
    -\frac{1}{2}\sum_{i,j}\mathcal{R}_{\min}g_{0}(f_{\ast}\bar{\omega}_{i},\omega_{j})(\bar{\omega}_{i}\otimes \omega_{j})&\geq& -\frac{2n(2n-1)|\mathcal{R}_{\min}|}{4};\label{eqn-ineqq}\\
    -\frac{1}{2}\sum_{i,j}g_{0}((\mathcal{R}-\mathcal{R}_{\min})(f_{\ast}\bar{\omega}_{i}),\omega_{j})(\bar{\omega}_{i}\otimes \omega_{j})&\geq&-\frac{\operatorname{Sc}_{g_{0}}-2n(2n-1)\mathcal{R}_{\min}}{4}.\label{eqn-ine}
\end{eqnarray}
Now we prove inequality (\ref{eqn-ine}). Since the operator $\mathcal{R}-\mathcal{R}_{\min}$ is nonnegative, there exists a self-adjoint operator $L$ on $\wedge^2 TM$ such that 
\[g_{0}((\mathcal{R}-\mathcal{R}_{\min})\omega_{i},\omega_{j})=g_{0}(L\omega_{i},L\omega_{j}).\]
Let
\[\bar{L}\bar{\omega}_{k}:=\sum_{i}g_{0}(L\omega_{k},f_{\ast}\bar{\omega}_{i})\bar{\omega}_{i}.\]
Then
\begin{eqnarray*}
   && -\sum_{i,j}g_{0}((\mathcal{R}-\mathcal{R}_{\min})f_{\ast}\bar{\omega}_{j},\omega_{i})[\bar{\omega}_{j}\otimes \omega_{i}]\\
   &=&-\sum_{i,j,k}g_{0}(Lf_{\ast}\bar{\omega}_{j},\omega_{k})g_{0}(L\omega_{i},\omega_{k})\bar{\omega}_{j}\otimes \omega_{i}\\
   &=&-\sum_{k}\bar{L}\bar{\omega}_{k}\otimes L\omega_{k}\\
   &=&\frac{1}{2}\sum_{k}\left[-(\bar{L}\bar{\omega}_{k}\otimes 1+1\otimes L\omega_{k})^2+(\bar{L}\bar{\omega}_{k})^2\otimes 1+1\otimes (L\omega_{k})^2\right]\\
   &\geq &\frac{1}{2}\sum_{k}\left[(\bar{L}\bar{\omega}_{k})^2\otimes 1+1\otimes (L\omega_{k})^2\right],
\end{eqnarray*}
since the Clifford multiplication by $2$-form is skew-symmetric and that squares of skew-symmetric endomorphisms are non-positive, so $-(\bar{L}\bar{\omega}_{k}\otimes 1+1\otimes L\omega_{k})^2$ is a non-negative endomorphism. We also have
\begin{eqnarray*}
  \sum_{k}(\bar{L}\bar{\omega}_{k})^2 &=& \sum_{i,j,k}g_{0}(L\omega_{k},f_{\ast}\bar{\omega}_{i})g_{0}(L\omega_{k},f_{\ast}\bar{\omega}_{j})\bar{\omega}_{i}\bar{\omega}_{j}\\
  &=&\sum_{i,j}g_{0}((\mathcal{R}-\mathcal{R}_{\min})f_{\ast}\bar{\omega}_{i},f_{\ast}\bar{\omega}_{j})\bar{\omega}_{i}\bar{\omega}_{j}.
\end{eqnarray*}
We can choose a local $g$-orthonormal basis $e_{1},\cdots,e_{2n}$ and a local $g_{0}$-orthonormal basic $\epsilon_{1},\cdots,\epsilon_{2n}$ such that
\[f_{\ast}e_{i}=\lambda_{i}\epsilon_{i}, \, \lambda_{i}\geq 0.\]
Then we have the orthonormal bases $e_{i}\wedge e_{j}$ of $\wedge^2TN$ and $\epsilon_{k}\wedge \epsilon_{l}$ of $\wedge^2TM$ and
\[f_{\ast}(e_{i}\wedge e_{j})=\lambda_{i}\lambda_{j}\epsilon_{i}\wedge \epsilon_{j}, \lambda_{i}\lambda_{j}\leq 1\]
for $1\leq i<j\leq 2n$, because $f$ is area nonincreasing.

\begin{eqnarray*}
   \sum_{k}(\bar{L}\bar{\omega}_{k})^2 &=& \sum_{i,j}g_{0}((\mathcal{R}-\mathcal{R}_{\min})f_{\ast}\bar{\omega}_{i},f_{\ast}\bar{\omega}_{j})\bar{\omega}_{i}\bar{\omega}_{j}\\
   &=&\sum_{i<j,k<l}g_{0}((\mathcal{R}-\mathcal{R}_{\min})f_{\ast}e_{i}\wedge f_{\ast}e_{j},f_{\ast}e_{k}\wedge f_{\ast}e_{l})e_{i}\wedge e_{j}\cdot e_{k}\wedge e_{l}\\
   &=&-\frac{1}{4}\sum_{i<j,k<l}\lambda_{i}\lambda_{j}\lambda_{k}\lambda_{l}g_{0}((\mathcal{R}-\mathcal{R}_{\min})\epsilon_{i}\wedge\epsilon_{j},\epsilon_{k}\wedge\epsilon_{l})e_{i}\cdot e_{j}\cdot e_{k}\cdot e_{l}\\
   &=&-\frac{1}{2}\sum_{i,j}\lambda^2_{i}\lambda_{j}^2 g_{0}((\mathcal{R}-\mathcal{R}_{\min})\epsilon_{i}\wedge\epsilon_{j},\epsilon_{i}\wedge\epsilon_{j})\text{ by equations (\ref{eqn-1-5})}\\
   &\geq &-\frac{1}{2}\sum_{i,j}g_{0}((\mathcal{R}-\mathcal{R}_{\min})\epsilon_{i}\wedge\epsilon_{j},\epsilon_{i}\wedge\epsilon_{j})\\
   &=&-\frac{\operatorname{Sc}_{g_{0}}-2n(2n-1)\mathcal{R}_{\min}}{2}.
\end{eqnarray*}
Similarly, we have
\[\sum_{k}(L\omega_{k})^2=-\frac{\operatorname{Sc}_{g_{0}}-2n(2n-1)\mathcal{R}_{\min}}{2}.\]
Hence, the inequality (\ref{eqn-ine}) follows.
For the proof of equality (\ref{eqn-ineqq}), we use the method of proof of Lemma 4.3 and Lemma 4.5 in \cite{ML1998}, then we have
\begin{eqnarray*}
   -\frac{1}{2}\mathcal{R}_{\min}\left\langle\sum_{i,j}g_{0}(f_{\ast}\bar{\omega}_{i},\omega_{j})(\bar{\omega}_{i}\otimes \omega_{j})\phi,\phi\right\rangle&\geq& -\left\langle\sum_{i\neq j}\frac{|\mathcal{R}_{\min}|\lambda_{i}\lambda_{j}}{4}\phi,\phi\right\rangle\\
   &\geq& -\frac{2n(2n-1)|\mathcal{R}_{\min}|}{4}\langle \phi,\phi\rangle.
\end{eqnarray*}
We have proved the claim in \eqref{eqn-ineqq} and \eqref{eqn-ine}.

For a harmonic spinor $\phi\in S$, by equation (\ref{eqn-Dirac}), we then have
\begin{multline*}
     \int_{N}\left(|\nabla \phi|^2+\left\langle \phi,\frac{\operatorname{Sc}_{g}-\operatorname{Sc}_{g_{0}}+2n(2n-1)\mathcal{R}_{\min}-2n(2n-1)|\mathcal{R}_{\min}|}{4}\phi\right\rangle\right)\leq 0.
\end{multline*}

Now if $\operatorname{Sc}_{g} (x) >\operatorname{Sc}_{g_{0}}(f(x))-2n(2n-1)\mathcal{R}_{\min}+2n(2n-1)|\mathcal{R}_{\min}|$, $\forall x \in N$, then $\phi=0$, so $\operatorname{Ind}(D_{E})=0$. This leads to a contradiction, since the Atiyah-Singer index theorem implies
\[\operatorname{Ind}(D_{E})=\operatorname{deg}(f)\chi(M)\neq 0,\]
for details, see the proof of Theorem 2.4 in \cite{GSSU2002}. Thus, the estimate in \eqref{ineqn-Sc} follows.

Moreover, if we assume $\operatorname{Sc}_{g}(x)\geq\operatorname{Sc}_{g_{0}}(f(x))-2n(2n-1)\mathcal{R}_{\min}+2n(2n-1)|\mathcal{R}_{\min}|$ for all $x \in N$, then $\nabla \phi=0$ and
\[\operatorname{Sc}_{g}(x)=\operatorname{Sc}_{g_{0}}(f(x))-2n(2n-1)\mathcal{R}_{\min}+2n(2n-1)|\mathcal{R}_{\min}|, \, \forall x \in N.\]
If $\mathcal{R}_{\min}\neq 0$, then all the inequalities in the proof of inequalities (\ref{eqn-ineqq}) and (\ref{eqn-ine}) become equalities. Tracking back the inequalities, we have $\lambda_{i}\lambda_{j}=1$ for all $i\neq j$. As a result, we have $\lambda_{i}=1$ for all $i=1,2,\cdots,2n$. Thus $f$ is a isometry and $\operatorname{Sc}_{g}(x)=\operatorname{Sc}_{g_{0}}(f(x)), \, \forall x \in N$. We also have 
\[-2n(2n-1)\mathcal{R}_{\min}+2n(2n-1)|\mathcal{R}_{\min}|=0\]
which implies $\mathcal{R}_{\min} \geq 0$.
\end{proof}


\section{Proof of Theorem \ref{Thm-manifold-boundary}}\label{sec-b}
In this section, we prove Theorem \ref{Thm-manifold-boundary}, basically following the argument in Lott \cite{JL2021} with some modification to remove the nonnegativity constraint on the curvature operator and the second fundamental form.
\begin{proof}
    We still consider the twisted spinor bundle $S=S(N)\otimes f^{\ast}S(M)$ on $N$.  Then for $\phi_{1},\phi_{2}\in S$, there is
    \begin{equation}\label{eqn-div}
        \int_{N}\langle D_E\phi_{1},\phi_{2}\rangle=\int_{N}\langle \phi_{1},D_E\phi_{2}\rangle+\int_{\partial N}\langle \nu_{N}\cdot\phi_{1},\phi_{2}\rangle
    \end{equation}
    for the outer unit normal vector $\nu_{N}$ on $\partial N$. Let $D^{\partial N}$ denote the Dirac operator on the boundary with respect to the restricted metric.
    For $\phi\in \Gamma(S)$ such that $D_E\phi=0$, by equations (\ref{eqn-Dirac}) and (\ref{eqn-div}), we have (equation (2.7) in \cite{JL2021})
     \begin{eqnarray}
        0&=&\int_{N}\langle \nabla\phi,\nabla\phi\rangle+\int_{\partial N}\langle \phi,D^{\partial N}\phi\rangle+\frac{1}{4}\int_{N}\operatorname{Sc}_{g}|\phi|^2 \nonumber\\
        &&+\int_{N}\langle \frac{1}{2}\sum_{j,k}(e_{j}\cdot e_{k}\cdot\sigma)\otimes R^{E}_{e_{j},e_{k}}(v),\phi\rangle \label{eqn-integral-LWB}\\
        &&+\frac{1}{2}\int_{\partial N}H_{N}\langle \phi,\phi\rangle-\frac{1}{2}\sum_{\alpha,i=1}^{2n-1}\int_{\partial N}A^{M}_{\alpha i}\langle \phi,\nu_{N}e_{i}\otimes\nu_{M}\epsilon_{\alpha}\phi\rangle, \nonumber
    \end{eqnarray}
    where $\{e_{i}\}_{i=1}^{2n}$ is a local orthonormal basis at point $x\in \partial N$ and $\{\epsilon_{\alpha}\}_{\alpha=1}^{2n}$ is a local orthonormal basic at point $f(x)\in \partial M$ such that $e_{2n}(x)=\nu_{N}(x)$ and $\epsilon_{2n}(f(x))=\nu_{M}(f(x))$, $A^M$ is the second fundamental form of $\partial M$ and $A^M_{\alpha i} = A^M((\partial f)_*(e_i), \epsilon_{\alpha})$.

    The term $\int_{N}\langle \frac{1}{2}\sum_{j,k}(e_{j}\cdot e_{k}\cdot\sigma)\otimes R^{E}_{e_{j},e_{k}}(v),\phi\rangle$ can be estimated as in the proof of Theorem \ref{Thm-main}, and we have
    \begin{eqnarray}\label{inequ-R}\nonumber
        \int_{N}\langle \frac{1}{2}\sum_{j,k}(e_{j}\cdot e_{k}\cdot\sigma)\otimes R^{E}_{e_{j},e_{k}}(v),\phi\rangle&\geq& -\frac{1}{4}\int_{N}(\operatorname{Sc}_{g_{0}}-2n(2n-1)\mathcal{R}_{\min})\langle \phi,\phi\rangle\\
        &&-\frac{2n(2n-1)|\mathcal{R}_{\min}|}{4}\int_{N}\langle \phi,\phi\rangle
    \end{eqnarray}
    For the term $\frac{1}{2}\int_{\partial N}A^{M}_{\alpha i}\langle \phi,\nu_{N}e_{i}\otimes\nu_{M}\epsilon_{\alpha}\phi\rangle$, we have
    \begin{eqnarray*}
        &&-\frac{1}{2}\int_{\partial N}\langle \phi,A^{M}_{\alpha i}e_{2n}e_{i}\otimes\epsilon_{2n}\epsilon_{\alpha}\phi\rangle\\
        &=&-\frac{1}{2}\int_{\partial N}\langle \phi,[A^{M}_{\alpha i}-A_{\min}^{M}(g_{0})_{\alpha i}+A_{\min}^{M}(g_{0})_{\alpha i}]e_{2n}e_{i}\otimes\epsilon_{2n}\epsilon_{\alpha}\phi\rangle\\
        &=&-\frac{1}{2}\int_{\partial N}\langle \phi,[A^{M}_{\alpha i}-A_{\min}^{M}(g_{0})_{\alpha i}]e_{2n}e_{i}\otimes\epsilon_{2n}\epsilon_{\alpha}\phi\rangle\\
        &&-\frac{1}{2}\int_{\partial N}\langle \phi,A_{\min}^{M}(g_{0})_{\alpha i}e_{2n}e_{i}\otimes\epsilon_{2n}\epsilon_{\alpha}\phi\rangle.
    \end{eqnarray*}
    Now we claim that
    \begin{eqnarray}
        & & -\frac{1}{2}\sum_{\alpha,i=1}^{2n-1}\int_{\partial N}\langle \phi,[A^{M}_{\alpha i}-A_{\min}^{M}(g_{0})_{\alpha i}]e_{2n}e_{i}\otimes\epsilon_{2n}\epsilon_{\alpha}\phi\rangle \nonumber \\
        &  \geq & -\frac{1}{2}\int_{\partial N}|\partial f|[H_{M}-(2n-1)A_{\min}^{M}]\langle \phi,\phi\rangle;\label{ineq-A}
    \end{eqnarray}
    \begin{equation}\label{ineq-A1}
    -\frac{1}{2}\sum_{\alpha,i=1}^{2n-1}\int_{\partial N}\langle \phi,A_{\min}^{M}(g_{0})_{\alpha i}e_{2n}e_{i}\otimes\epsilon_{2n}\epsilon_{\alpha}\phi\rangle\geq -\frac{1}{2}|\partial f|A_{\min}^{M}(2n-1)\int_{\partial N}\langle \phi,\phi\rangle.
    \end{equation}
    We first check inequality (\ref{ineq-A}). Following the argument in Lemma 2.1 in \cite{JL2021} and also Lemma 2.3 in \cite{WX2023}, since the operator $A^{M}-A_{\min}^{M}g_{0}$ is non-negative, there exists a self-adjoint operator $L\in \operatorname{End}(T\partial M)$ such that $(A^{M}-A_{\min}^{M}g_{0})(\cdot,\cdot)=g_{0}(L\cdot,L\cdot)$. Let 
    \[\bar{L}\epsilon_{\beta}=\sum_{i=1}^{2n-1}g_{0}(L\epsilon_{\beta},f_{\ast}e_{i})e_{i}\]
    Then, for $\theta>0$,
    \begin{eqnarray*}
        &&-\frac{1}{2}\sum_{i,\alpha=1}^{2n-1}[A^{M}_{\alpha i}-A_{\min}^{M}(g_{0})_{\alpha i}]e_{2n}e_{i}\otimes\epsilon_{2n}\epsilon_{\alpha}\\
        &=&-\frac{1}{2}\sum_{i,\alpha}g_{0}(L(f_{\ast}e_{i}),L\epsilon_{\alpha})e_{2n}e_{i}\otimes\epsilon_{2n}\epsilon_{\alpha}\\
        &=&-\frac{1}{2}\sum_{i,\alpha,\beta}g_{0}(L(f_{\ast}e_{i}),\epsilon_{\beta})g_{0}(\epsilon_{\beta},L\epsilon_{\alpha})e_{2n}e_{i}\otimes\epsilon_{2n}\epsilon_{\alpha}\\
        &=&-\frac{1}{2}\sum_{\beta}e_{2n}\bar{L}\epsilon_{\beta}\otimes\epsilon_{2n}L\epsilon_{\beta}\\
        &=&\sum_{\beta}\frac{1}{4}\left[\theta^{-2}(e_{2n}\bar{L}\epsilon_{\beta})^2\otimes 1+1\otimes \theta^2(\epsilon_{2n}L\epsilon_{\beta})^2-(\theta^{-1}e_{2n}L\epsilon_{\beta}\otimes 1+1\otimes \theta\epsilon_{2n}L\epsilon_{\beta})^2 \right]\\
        &\geq&\sum_{\beta}\frac{1}{4}\left[\theta^{-2}(e_{2n}\bar{L}\epsilon_{\beta})^2\otimes 1+1\otimes \theta^{2}(\epsilon_{2n}L\epsilon_{\beta})^2\right],
    \end{eqnarray*}
    since $(\theta^{-1}e_{2n}\bar{L}\epsilon_{\beta}\otimes 1+1\otimes \theta\epsilon_{2n}L\epsilon_{\beta})$ is skew-symmetric. We also have
    \begin{eqnarray*}
        (\epsilon_{2n}L\epsilon_{\beta})^2&=&\sum_{\beta,\alpha,\gamma}g_{0}(L\epsilon_{\beta},\epsilon_{\alpha})g_{M}(L\epsilon_{\beta},\epsilon_{\gamma})\epsilon_{2n}\epsilon_{\alpha}\epsilon_{2n}\epsilon_{\gamma}\\
        &=&\sum_{\beta,\alpha}g_{0}(L\epsilon_{\beta},\epsilon_{\alpha})g_{0}(L\epsilon_{\beta},\epsilon_{\alpha})\epsilon_{2n}\epsilon_{\alpha}\epsilon_{2n}\epsilon_{\alpha}\\
        &&+\sum_{\beta,\alpha\neq\gamma}g_{0}(L\epsilon_{\beta},\epsilon_{\alpha})g_{0}(L\epsilon_{\beta},\epsilon_{\gamma})\epsilon_{2n}\epsilon_{\alpha}\epsilon_{2n}\epsilon_{\gamma}\\
        &=&-\sum_{\beta}g_{0}(L\epsilon_{\beta},L\epsilon_{\beta})+\sum_{\alpha\neq\gamma}g_{0}(L\epsilon_{\gamma},L\epsilon_{\alpha})\epsilon_{\alpha}\epsilon_{\gamma}\\
        &=&-\sum_{\beta}(A^{M}-A_{\min}^{M}g_{0})(\epsilon_{\beta},\epsilon_{\beta})\\
        &=&-(H_{M}-(2n-1)A_{\min}^{M}).
    \end{eqnarray*}
    Moreover, 
    \begin{eqnarray*}
        \sum_{\beta}(e_{2n}\bar{L}\epsilon_{\beta})^2&=&\sum_{\beta,i,j=1}^{2n-1}g_{0}(L\epsilon_{\beta},f_{\ast}e_{i})g_{0}(L\epsilon_{\beta},f_{\ast}e_{j})e_{2n}e_{i}e_{2n}e_{j}\\
        &=&\sum_{i,j=1}^{2n-1}g_{0}(Lf_{\ast}e_{j},Lf_{\ast}e_{i})e_{2n}e_{i}e_{2n}e_{j}\\
        &=&-\sum_{i}^{2n-1}g_{0}(Lf_{\ast}e_{i},Lf_{\ast}e_{i})+\sum_{i\neq j}g_{0}(Lf_{\ast}e_{j},Lf_{\ast}e_{i})e_{i}e_{j}\\
        &=&-\sum_{i}(A^{M}-A_{\min}^{M}g_{0})(f_{\ast}e_{i},f_{\ast}e_{i})\\
        &\geq&-|\partial f|^2(H_{M}-(2n-1)A_{\min}^{M})
    \end{eqnarray*}
    where 
    \[|\partial f|^2:=\sup_{0\neq v\in T_{p}\partial N}\frac{g_{0}(\partial f_{\ast}v,\partial f_{\ast}v)}{g(v,v)},  p\in \partial N.\]
    Choosing $\theta=\sqrt{|\partial f|}$ for $|\partial f|\neq 0$, we have
      \begin{eqnarray*}
        &&-\frac{1}{2}\sum_{\alpha,i=1}^{2n-1}\int_{\partial N}\langle \phi,[A^{M}_{\alpha i}-A_{\min}^{M}(g_{0})_{\alpha i}]e_{2n}e_{i}\otimes\epsilon_{2n}\epsilon_{\alpha}\phi\rangle\\
        &\geq& -\frac{1}{2}\int_{\partial N}|\partial f|[H_{M}-(2n-1)A^{M}_{\min}]\langle \phi,\phi\rangle.  
      \end{eqnarray*}
    Now we check the inequality (\ref{ineq-A1}). We can compute in local coordinate as in Lemma 4.5 in \cite{ML1998}. We can choose the local $\{e_{i}\}_{i=1}^{2n-1}$ and $\{\epsilon_{\alpha}\}_{\alpha=1}^{2n-1}$ such that $(\partial f)_{\ast}e_{i}=\mu_{i}\epsilon_{i}$ for $i=1,\cdots,2n-1$. Then we have
    \begin{eqnarray*}
      -\frac{1}{2}\sum_{\alpha,i=1}^{2n-1}\langle A_{\min}^{M}(g_{0})_{\alpha i}e_{2n}e_{i}\otimes\epsilon_{2n}\epsilon_{\alpha}\phi,\phi\rangle&\geq& -\frac{1}{2}|A^{M}_{\min}|\left(\sum_{i=1}^{2n-1}\mu_{i}\right)\langle \phi,\phi\rangle  \\
     &\geq&  -\frac{1}{2}(2n-1)|A^{M}_{\min}||\partial f|\langle \phi,\phi\rangle.  
    \end{eqnarray*}
    Integrating over $\partial N$, the  inequality (\ref{ineq-A1}) follows. 
    
    Then by plugging (\ref{inequ-R}), (\ref{ineq-A}) and (\ref{ineq-A1}) into (\ref{eqn-integral-LWB}), we obtain 
    \begin{eqnarray*}
        0&\geq&\int_{N}\langle \nabla\phi,\nabla\phi\rangle+\int_{\partial N}\langle \phi,D^{\partial N}\phi\rangle+\frac{1}{4}\int_{N}\operatorname{Sc}_{g}|\phi|^2\\
        &&-\frac{1}{4}\int_{N}(\operatorname{Sc}_{g_{0}}-2n(2n-1)\mathcal{R}_{\min})\langle \phi,\phi\rangle-\frac{2n(2n-1)|\mathcal{R}_{\min}|}{4}\int_{N}\langle \phi,\phi\rangle\\
        &&+\frac{1}{2}\int_{\partial N}H_{N}\langle \phi,\phi\rangle-\frac{1}{2}\int_{\partial N}|\partial f|[f^{\ast}H_{M}-(2n-1)A_{\min}^{M}]\langle \phi,\phi\rangle\\
        &&-\frac{1}{2}|\partial f| (2n-1) |A_{\min}^{M}| \int_{\partial N}\langle \phi,\phi\rangle.
    \end{eqnarray*}

    We impose the boundary condition as in section 2.2 in \cite{JL2021} or section 3.1 in \cite{WX2023} such that $\int_{\partial N}\langle \phi,D^{\partial N}\phi\rangle\geq 0$ and $\operatorname{Ind}(D_{E})\neq 0$, for details, we refer to \cite{JL2021} and \cite{WX2023}. 
    Then because we assume that $\partial f$ is distance nonincreasing, $|\partial f|\leq 1$. Note that $f^* H_M - (2n-1)A^{M}_{\min} \geq 0$. If 
    \[\left[\operatorname{Sc}_{g}(x)-\operatorname{Sc}_{g_{0}}(f(x))\right]> -2n(2n-1)\mathcal{R}_{\min}+2n(2n-1)|\mathcal{R}_{\min}|\]   and
    \[
        \left[H_{ N}(y)-H_{ M}(\partial f(y))\right]> -(2n-1)A^{M}_{\min}+(2n-1)|A^{M}_{\min}|.
    \]
    Then we have $\phi=0$ and $\operatorname{Ind}(D_{E})= 0$ which leads to a contradiction. Thus, either \eqref{eqn-scal-bdry} or \eqref{eqn-mean-curv-bdry} must hold.

    Finally, if  \[\left[\operatorname{Sc}_{g}(x)-\operatorname{Sc}_{g_{0}}(f(x))\right]\geq -2n(2n-1)\mathcal{R}_{\min}+2n(2n-1)|\mathcal{R}_{\min}|, \, \forall x \in N\]   and
    \[
        \left[H_{ N}(y)-H_{ M}(\partial f(y))\right]\geq -(2n-1)A^{M}_{\min}+(2n-1)|A^{M}_{\min}|, \, \forall y \in \partial N.
    \]
    Then \[\left[\operatorname{Sc}_{g}(x)-\operatorname{Sc}_{g_{0}}(f(x))\right]= -2n(2n-1)\mathcal{R}_{\min}+2n(2n-1)|\mathcal{R}_{\min}|\, \forall x \in N\]
    and
    \[
        \left[H_{ N}(y)-H_{ M}(\partial f(y))\right]= -(2n-1)A^{M}_{\min}+(2n-1)|A^{M}_{\min}|, \, \forall y \in \partial N
    \]
    and $\phi$ is a parallel twisted spinor. All the inequalities become equalities in the proof of Theorem \ref{Thm-manifold-boundary}. By tracking those inequalities, we obtain that
    \[\operatorname{Sc}_{g}(x)=\operatorname{Sc}_{g_{0}}(f(x)), \forall x \in N, \quad \mathcal{R}_{\min}\geq 0\]
    and
    \[
    \begin{cases}
        H_{ N}(y)=H_{ M}(\partial f(y)), \forall y \in \partial N, &\quad \text{if } A^{M}_{\min}= 0\\
        \partial f:\partial N\to\partial M \text{ is local isometric }, &\quad \text{if } A^{M}_{\min}\neq 0
    \end{cases}
    \]
\end{proof}
\begin{remark}
{\rm
    In Corollary \ref{cor-Euclidean domain}, since the metric $g_{E}$ is flat, the curvature term $R^{E}$ in \eqref{eqn-Dirac} is zero. Therefore, we do not need to handle that term.
}
\end{remark}

\bibliographystyle{plain}
\bibliography{main.bib}

\end{document}